\documentclass[12pt]{article}
\usepackage[T2A]{fontenc}
\usepackage[utf8]{inputenc}
\usepackage{amssymb,amsmath,amsfonts,amsthm,amscd,latexsym,indentfirst,verbatim}
\def\Span{\mathrm{Span}}
\def\Aut{\mathrm{Aut}}

\begin{document}

\begin{center}
{\Large Nonunital decompositions of the matrix algebra of order three}

V. Gubarev
\end{center}

\begin{abstract}
All decompositions of $M_3(\mathbb{C})$ into a direct vector-space sum of two
subalgebras such that none of the subalgebras contains the identity matrix are classified.
Thus, the classification of all decompositions of $M_3(\mathbb{C})$ 
into a direct vector-space sum of two subalgebras as well as 
description of Rota---Baxter operators of nonzero weight on $M_3(\mathbb{C})$ is finished.

\medskip
{\it Keywords}: decomposition of algebra, sum of rings, matrix algebra.

MSC: 16S50
\end{abstract}

Let $R$ be a ring, suppose that $R = R_1 + R_2$, where $R_1,R_2$ are subrings (not necessarily ideals) of $R$. 
In this situation, we say that $R$ decomposes into a sum of $R_1,R_2$.
If $R_1\cap R_2 = (0)$, then we call such decomposition as a direct one.
The study of decompositions of associative rings and algebras started in 1963, when O.H.~Kegel proved~\cite{Kegel1963} that an associative algebra 
decomposed into a sum of two nilpotent subalgebras is itself nilpotent.

In 1995, K.I. Beidar and A.V. Mikhalev asked~\cite{BeidarMikhalev95}, if a sum of two PI-algebras is again a~PI-algebra.
In 2017, this problem was positively solved by M.~K\k{e}pczyk~\cite{Kepczyk2017}.
The famous K\"{o}the problem (If a ring $R$ has no nonzero nil ideals,
does it follow that $R$~has no nonzero nil one-sided ideals?)
is equivalent to a problem concerned decompositions~\cite{Ferrero1989}.
In~\cite{Kosan}, all decompositions of 2- and 3-dimensional associative algebras were decribed.

It is natural to study decompositions involved matrix algebras.
In~1999, Y.A.~Bahturin and O.H.~Kegel~\cite{BahturinKegel} described all algebras
decomposed as a sum of two matrix algebras 
(it is Problem~3 from~\cite{Kosan} posted twenty years later in 2019). 

In~\cite{GonGub}, all direct decompositions of $M_2(\mathbb{C})$
were classified. 
In the current work we finish a classification
of direct decompositions of $M_3(\mathbb{C})$.
In~\cite{Gub}, all direct decompositions of $M_3(\mathbb{C})$
such that one of the subalgebras contains the identity matrix were classified (71 cases, some of them involve one or two parameters).
Thus, it remains to describe direct decompositions of $M_3(\mathbb{C})$
such that none of the two subalgebras contains the identity matrix (we call them as nonunital ones). 
Joint with the articles~\cite{GonGub,Gub}, we complete the description of all Rota---Baxter operators (algebraic analogue of integral operator~\cite{Guo2011}) of nonzero weight on $M_3(\mathbb{C})$.

We split the problem of classification of a nonunital decomposition
$M_3(\mathbb{C}) = S\oplus M$, where $\dim M>\dim S$, into the following cases:
$\dim M = 6$ or $\dim M = 5$.
In the first case, there is the only possibility for $M$~\cite{Unital}: 
it is the subalgebra of matrices with zero first column. 
In the second one, due to~\cite{Gub}, we have either 
$M = \Span\{e_{11},e_{12},e_{13},e_{22},e_{23}\}$ or 
$M = \Span\{e_{11},e_{12},e_{13},e_{23},e_{23}\}$.
We actively apply automorphisms of $M_3(\mathbb{C})$, which preserve the subalgebra~$M$, it helps us to avoid computational difficulties.

In Theorem~1, we get three decompositions, when $\dim M = 6$.
In Theorem~2, we obtain nine more decompositions, here $\dim M = 5$.
Note that all twelve cases do note involve parameters.
Therefore, due to~\cite{GonGub,Gub}, we have 83 nontrivial decompositions of~$M_3(\mathbb{C})$
and 119 nontrivial Rota--Baxter operators on~$M_3(\mathbb{C})$.

In what follows, we will apply an automorphism $\Theta_{12}$ of $M_3(\mathbb{C})$,
acting as follows, $\Theta_{12}(X) = T^{-1}XT$ for $T = e_{12}+e_{21}+e_{33}$.
Analogously, we define $\Theta_{13},\Theta_{23}\in\Aut(M_3(\mathbb{C}))$.

\section{(6,3)-decompositions}

Let $M$ be a six-dimensional nonunital subalgebra of $M_3(\mathbb{C})$.
Then $M$ consists only of degenerate matrices, and then up to transpose~$M$ is 
isomorphic to the subalgebra having all zero elements in the first column~\cite{Meshulam}.
Thus, we fix six-dimensional subalgebra $M = \Span\{e_{12},e_{13},e_{22},e_{23},e_{32},e_{33}\}$.

{\bf Lemma 1}~\cite{Gub}.
Let $\varphi$ be an automorphism of $M_3(\mathbb{C})$ preserving 
the subalgebra $M = \Span\{e_{12},e_{13},e_{22},e_{23},e_{32},e_{33}\}$.
Then $\varphi$ acts as follows,
\begin{gather}
e_{11}\to \begin{pmatrix}
1 & \beta & \gamma \\
0 & 0 & 0 \\
0 & 0 & 0 \\
\end{pmatrix},\quad
e_{12}\to \begin{pmatrix}
0 & \kappa & \lambda \\
0 & 0 & 0 \\
0 & 0 & 0 \\
\end{pmatrix},\quad
e_{13}\to \begin{pmatrix}
0 & \mu & \nu \\
0 & 0 & 0 \\
0 & 0 & 0 \\
\end{pmatrix}, \nonumber \\ 
e_{22}\to \frac{1}{\Delta}\begin{pmatrix}
0 & \kappa(\gamma \mu-\beta \nu) & \lambda(\gamma \mu-\beta \nu) \\
0 & \kappa \nu & \lambda \nu \\
0 & -\kappa \mu & -\lambda \mu \\
\end{pmatrix},\
e_{23}\to \frac{1}{\Delta}\begin{pmatrix}
0 & \mu(\gamma \mu-\beta \nu) & \nu(\gamma \mu-\beta \nu) \\
0 & \mu\nu & \nu^2 \\
0 & -\mu^2 & -\mu\nu \\
\end{pmatrix}, \nonumber \\
e_{32}\to \frac{1}{\Delta}\begin{pmatrix}
0 & \kappa(\beta \lambda-\gamma \kappa) & \lambda(\beta \lambda-\gamma \kappa) \\
0 & -\kappa\lambda & -\lambda^2 \\
0 & \kappa^2 & \kappa\lambda \\
\end{pmatrix},\
e_{33}\to \frac{1}{\Delta}\begin{pmatrix}
0 & \mu(\beta \lambda-\gamma \kappa) & \nu(\beta \lambda-\gamma \kappa) \\
0 & -\lambda\mu & -\lambda\nu \\
0 & \kappa\mu & \kappa\nu \\
\end{pmatrix}, \nonumber \\
e_{21}\to \frac{1}{\Delta}\begin{pmatrix}
\gamma \mu - \beta \nu & \beta(\gamma \mu -\beta \nu) & \gamma(\gamma \mu-\beta \nu) \\
\nu & \beta \nu & \gamma \nu \\
-\mu & -\beta \mu & -\gamma \mu  \\
\end{pmatrix}, \label{AutPreserv} \\
e_{31}\to \frac{1}{\Delta}\begin{pmatrix}
\beta \lambda-\gamma \kappa & \beta(\beta \lambda-\gamma \kappa) 
 & \gamma(\beta \lambda-\gamma \kappa) \\
-\lambda & -\beta \lambda & -\gamma \lambda \\
\kappa & \beta \kappa & \gamma \kappa \\
\end{pmatrix}, \nonumber
\end{gather}
where $\Delta = \kappa\nu-\lambda\mu\neq0$.

{\bf Theorem 1}.
Every direct decomposition of $M_3(\mathbb{C})$
with two subalgebras of the dimensions~3 and~6 
not containing the identity matrix, up to transpose
and up to action of $\Aut(M_3(\mathbb{C}))$ is isomorphic to
one of the following cases:

(A1) $\Span\{e_{11},e_{21},e_{31}\}
 \oplus \Span\{e_{12},e_{13},e_{22},e_{23},e_{32},e_{33}\}$,

(A2) $\Span\{e_{11}+e_{22},e_{21},e_{31}\}
 \oplus \Span\{e_{12},e_{13},e_{22},e_{23},e_{32},e_{33}\}$,

(A3) $\Span\{e_{11}+e_{22},e_{21}+e_{22},e_{31}\}
 \oplus \Span\{e_{12},e_{13},e_{22},e_{23},e_{32},e_{33}\}$.

Moreover, all cases (A1)--(A3) lie in different orbits under action 
automorphisms or antiautomorphisms of $M_3(\mathbb{C})$.

{\sc Proof}. 
We study a decomposition $M_3(\mathbb{C}) = S\oplus M$, and $S$ has a~basis
$$
v_1 = \begin{pmatrix}
0 & a & b \\
1 & c & d \\
0 & e & f \\
\end{pmatrix}, \quad
v_2 = \begin{pmatrix}
0 & r & s \\
0 & t & u \\
1 & x & y \\
\end{pmatrix}, \quad
v_3 = \begin{pmatrix}
1 & k & l \\
0 & m & n \\
0 & p & q \\
\end{pmatrix}.
$$
Note that $v_3^2 = v_3$, since $S$ is a~subalgebra. 
Hence, its submatrix 
$H = \begin{pmatrix}
m & n \\
p & q \\
\end{pmatrix}$ is idempotent.
We may conjugate $H$ with appropriate nondegenerate matrix~$T$
such that $J = T^{-1}HT$ is in the Jordan form. 
We have the following variants for $J\in M_2(\mathbb{C})$:

1) $J = 0$,

2) $J = \begin{pmatrix}
1 & 0 \\
0 & 0 \\
\end{pmatrix}$,

3) $J = \begin{pmatrix}
0 & 0 \\
0 & 1 \\
\end{pmatrix}$,

4) $J = E$.

Let us conjugate $v_3$ with the block-diagonal matrix
$\begin{pmatrix}
1 & 0 \\
0 & T \\
\end{pmatrix}$,
such conjugation surely preserves~$M$.
In the case 4) involving the equality $v_3^2 = v_3$, we get $v_3 = E$, a~contradiction.

The case 3) is conjugate to 2) under $\Theta_{23}$.
Thus, it remains to consider cases~1) and~2).

{\sc Case 1}: $J = 0$.
Hence, $v_3 = e_{11} + ke_{12} + le_{13}$.
Note that the automorphism~$\varphi$ defined by~\eqref{AutPreserv}
with $\beta = k$ and $\gamma = l$ maps $e_{11}$ to $v_3$.
So, $\varphi^{-1}$ maps $v_3$ to $e_{11}$, and we may assume that $v_3 = e_{11}$.
Since $v_1e_{11},v_2e_{11}\in S$, we get the decomposition (A1).

{\sc Case 2}: $J = \begin{pmatrix}
1 & 0 \\
0 & 0 \\
\end{pmatrix}$.
The condition $v_3^2 = v_3$ implies that
$v_3 = \begin{pmatrix}
1 & 0 & l \\
0 & 1 & 0 \\
0 & 0 & 0 \\
\end{pmatrix}$.
Again, we apply~$\varphi$ defined by~\eqref{AutPreserv}
with $\beta = \mu = \lambda = 0$ and $\gamma = l$,
then $\varphi(e_{11}+e_{22}) = v_3$.
We may assume that $v_3 = e_{11}+e_{22}$.
Since $v_iv_3,v_3v_i\in S$ for $i=1,2$, we conclude that
$$
v_1 = \begin{pmatrix}
0 & a & 0 \\
1 & c & 0 \\
0 & 0 & 0 \\
\end{pmatrix}, \quad
v_2 = \begin{pmatrix}
0 & 0 & 0 \\
0 & 0 & 0 \\
1 & x & 0 \\
\end{pmatrix}.
$$

Taking $\varphi$~\eqref{AutPreserv}
with $\lambda = \mu = \gamma = 0$, $\nu = 1$, and $\beta = x$,
then $\varphi(v_3) = v_3$
and $\varphi(e_{31}) = v_2$.
Therefore, we may assume that $v_2 = e_{31}$.
From $v_2v_1\in S$, we get $a = 0$,
If $c = 0$, then it is exactly the decomposition~(A2).
Otherwise, consider $\varphi$~\eqref{AutPreserv}
with $\lambda = \mu = \beta = \gamma = 0$, $\kappa = c$, $\nu = 1$,
then $\varphi(v_2) = v_2$, $\varphi(v_3) = v_3$,
and $\varphi(e_{21}+e_{22}) = (1/c)v_1$.
Hence, we arrive at the decomposition~(A3).

Decompositions~(A1) and~(A2) as well as (A1) and (A3) lie in different orbits 
under action (anti)automorphisms of $M_3(\mathbb{C})$ 
preserving~$M$,
since the 3-dimensional subalgebra from~(A2) and (A3) 
but not from~(A1) contains an idempotent of rank~2.
Finally, decompositions~(A2) and (A3) lie in different orbits too.
Indeed, the dimensions of the semisimple parts of the subalgebra~$S$ 
from (A2) and (A3) do not equal. 
\hfill $\square$

\section{(5,4)-decompositions}

Let $M$ be a nonunital 5-dimensional algebra,
by \cite[Lemma~5]{Gub}, 
we may assume up to transpose and action of $\Aut(M_3(\mathbb{C}))$
that either $M = \Span\{e_{11},e_{12},e_{13},e_{22},e_{23}\}$ or
$M = \Span\{e_{11},e_{12},e_{13},e_{23},e_{33}\}$.
Then the group $\Aut(M_3(\mathbb{C}))$ preserving~$M$
coincides with the group of automorphisms of $\Aut(M_3(\mathbb{C}))$
preserving the subalgebra of upper-triangular matrices.
Thus, an automorphism $\psi\in\Aut(M_3(\mathbb{C}))$ preserving~$M$
has the form~\eqref{AutPreserv} considered with $\mu = 0$:
\begin{equation}\label{AutU}
\begin{gathered}
e_{11}\to \begin{pmatrix}
1 & \beta & \gamma \\
0 & 0 & 0 \\
0 & 0 & 0 \\
\end{pmatrix},\quad
e_{12}\to \begin{pmatrix}
0 & \delta & \varepsilon \\
0 & 0 & 0 \\
0 & 0 & 0 \\
\end{pmatrix},\quad
e_{13}\to \begin{pmatrix}
0 & 0 & \alpha \\
0 & 0 & 0 \\
0 & 0 & 0 \\
\end{pmatrix}, \\
e_{21}\to
\frac{1}{\delta}\begin{pmatrix}
-\beta & -\beta^2 & -\beta\gamma \\
1 & \beta & \gamma \\
0 & 0 & 0 \\
\end{pmatrix},\quad
e_{22}\to
\begin{pmatrix}
0 & -\beta & -\beta\varepsilon/\delta \\
0 & 1 & \varepsilon/\delta \\
0 & 0 & 0 \\
\end{pmatrix}, \\
e_{23}\to
\frac{1}{\delta}\begin{pmatrix}
0 & 0 & -\alpha\beta \\
0 & 0 & \alpha \\
0 & 0 & 0 \\
\end{pmatrix},\quad
e_{31}\to
\frac{1}{\alpha\delta}
\begin{pmatrix}
\beta\varepsilon-\gamma\delta & \beta(\beta\varepsilon-\gamma\delta)
 & \gamma(\beta\varepsilon-\gamma\delta) \\
-\varepsilon & -\beta\varepsilon & -\gamma\varepsilon \\
\delta & \beta\delta & \gamma\delta \\
\end{pmatrix}, \\
e_{32}\to
\frac{1}{\alpha}\begin{pmatrix}
0 & \beta\varepsilon-\gamma\delta & \varepsilon(\beta\varepsilon-\gamma\delta)/\delta \\
0 & -\varepsilon & -\varepsilon^2/\delta \\
0 & \delta & \varepsilon \\
\end{pmatrix}, \quad
e_{33}\to
\begin{pmatrix}
0 & 0 & \beta\varepsilon/\delta -\gamma \\
0 & 0 & -\varepsilon/\delta \\
0 & 0 & 1 \\
\end{pmatrix},
\end{gathered}
\end{equation}
where $\alpha,\delta \neq0$.

{\bf Theorem 2}.
Every direct decomposition of $M_3(\mathbb{C})$
with two subalgebras of the dimensions~4 and~5 
not containing the identity matrix, up to transpose
and up to action of $\Aut(M_3(\mathbb{C}))$ is isomorphic to
one of the following cases:

(B1) $\Span\{e_{21},e_{31},e_{32},e_{33}\}\oplus 
\Span\{e_{11},e_{12},e_{13},e_{22},e_{23}\}$,

(B2) $\Span\{e_{11}+e_{21},e_{31},e_{32},e_{33}\}\oplus 
\Span\{e_{11},e_{12},e_{13},e_{22},e_{23}\}$,

(B3) $\Span\{e_{21},e_{31},e_{32},e_{22}+e_{33}\}\oplus 
\Span\{e_{11},e_{12},e_{13},e_{22},e_{23}\}$,

(B4) $\Span\{e_{21},e_{31},e_{32}+e_{23},e_{22}+e_{33}\}\oplus 
\Span\{e_{11},e_{12},e_{13},e_{22},e_{23}\}$;

(B5) $\Span\{e_{21},e_{31},e_{32},e_{11}+e_{33}\}\oplus 
\Span\{e_{11},e_{12},e_{13},e_{22},e_{23}\}$,

(B6) $\Span\{e_{21},e_{31},e_{32},e_{22}\}\oplus
\Span\{e_{11},e_{12},e_{13},e_{23},e_{33}\}$,

(B7) $\Span\{e_{21},e_{11}+e_{31},e_{12}+e_{32},e_{22}\}\oplus
\Span\{e_{11},e_{12},e_{13},e_{23},e_{33}\}$,

(B8) $\Span\{e_{21},e_{31},e_{32},e_{22}+e_{33}\}\oplus
\Span\{e_{11},e_{12},e_{13},e_{23},e_{33}\}$,

(B9) $\Span\{e_{21},e_{31},e_{32}+e_{23},e_{22}+e_{33}\}\oplus
\Span\{e_{11},e_{12},e_{13},e_{23},e_{33}\}$.

Moreover, all cases (B1)--(B9) lie in different orbits under action automorphisms or antiautomorphisms of $M_3(\mathbb{C})$.

{\sc Proof}. 
Let us start with the case $M = \Span\{e_{11},e_{12},e_{13},e_{22},e_{23}\}$.
Then the complement subalgebra $S$ contains a basis
$$
v_1 = \begin{pmatrix}
a & b & c \\
1 & d & e \\
0 & 0 & 0 \\
\end{pmatrix}, \quad
v_2 = \begin{pmatrix}
f & g & h \\
0 & i & j \\
1 & 0 & 0 \\
\end{pmatrix}, \quad
v_3 = \begin{pmatrix}
k & l & m \\
0 & n & p \\
0 & 1 & 0 \\
\end{pmatrix}, \quad
v_4 = \begin{pmatrix}
x & y & z \\
0 & t & u \\
0 & 0 & 1 \\
\end{pmatrix}.
$$

As in the proof of Theorem~1, we have $v_4^2 = v_4$.
It is known~\cite{Thijsse} that an upper-triangular matrix from
$M_3(\mathbb{C})$ is conjugate to its Jordan form with the help of 
some upper-triangular matrix, i.\,e. an automorphism preserving~$M$.
Thus, we have the following cases for the Jordan form~$J$ of $v_4$:

1) $J = e_{33}$,

2) $J = e_{22}+e_{33}$,

3) $J = e_{11}+e_{33}$.

{\sc Case 1}: $v_4 = e_{33}$.
From $v_4v_2 = v_2$ and $v_4v_3 = v_3$, we conclude that
$v_2 = e_{31}$ and $v_3 = e_{32}$.
Also, $v_1v_4 = 0$, hence, $c = e = 0$.
Since $v_1^2 = (a+d)v_1$, we get 
$v_1 = \begin{pmatrix}
a & ad & 0 \\
1 & d & 0 \\
0 & 0 & 0 \\
\end{pmatrix}$.
If $a+d = 0$, then we apply $\psi$ defined by~\eqref{AutU} with parameters
$\gamma = \varepsilon = 0$, $\beta = d$.
So, $\psi^{-1}$ gives us the decomposition~(B1).
If $a+d\neq0$, then consider $\psi$~\eqref{AutU} with parameters
$\gamma = \varepsilon = 0$, $\beta = d$ and $\delta = a+d\neq0$.
It is easy to check that $\psi(e_{11}+e_{21}) = v_1$, however,
$L(v_2,v_3,v_4)$ is $\psi$-invariant. It is the decomposition~(B2).

{\sc Case 2}: $v_4 = e_{22}+e_{33}$.
From the equalities 
$$
v_4v_1 = v_1, \quad v_1v_4 = 0, \quad
v_4v_2 = v_2, \quad v_2v_4 = 0, \quad
v_4v_3 = v_3v_4 = v_3,
$$
we get $v_1 = e_{21}$, $v_2 = e_{31}$ and
$v_4 = \begin{pmatrix}
0 & 0 & 0 \\
0 & n & p \\
0 & 1 & 0 \\
\end{pmatrix}$.
If $p+n^2/4 = 0$, we take $\psi$~\eqref{AutU} with 
$\beta = \gamma = 0$, $\varepsilon/\delta = -n/2$
and get the decomposition~(B3).
Otherwise, we take $\psi$~\eqref{AutU} with 
$\beta = \gamma = 0$, $\varepsilon = -n/2$,
$\delta = 1$, $\alpha = \sqrt{p+n^2/4}$, 
and get the decomposition~(B4).

{\sc Case 3}: $v_4 = e_{11}+e_{33}$.
From the equalities 
$$
v_4v_1 = 0, \quad 
v_1v_4 = v_1, \quad
v_4v_2 = v_2v_4 = v_2, \quad
v_4v_3 = v_3, \quad 
v_3v_4 = 0,
$$
we get 
$$
v_1 = \begin{pmatrix}
0 & 0 & 0 \\
1 & 0 & e \\
0 & 0 & 0 \\
\end{pmatrix}, \quad
v_2 = \begin{pmatrix}
f & 0 & h \\
0 & 0 & 0 \\
1 & 0 & 0 \\
\end{pmatrix}, \quad
v_3 = \begin{pmatrix}
0 & l & 0 \\
0 & 0 & 0 \\
0 & 1 & 0 \\
\end{pmatrix}.
$$
Further, the relations
$v_1v_3 = 0$ and $v_3v_1 = v_2 + ev_4$
imply 
$l+e = e + f - l = h - el = 0$.
Thus, $l = -e$, $f = -2e$, $h = -e^2$. 
We apply $\psi$~\eqref{AutU} with
$\beta = \varepsilon = 0$ and $\gamma = e$,
and get the decomposition~(B5) with the help of $\psi^{-1}$.

Now, we consider the case $M = \Span\{e_{11},e_{12},e_{13},e_{23},e_{33}\}$.
The complement subalgebra $S$ contains a basis
$$
v_1 = \begin{pmatrix}
a & b & c \\
1 & 0 & d \\
0 & 0 & e \\
\end{pmatrix}, \quad
v_2 = \begin{pmatrix}
f & g & h \\
0 & 0 & i \\
1 & 0 & j \\
\end{pmatrix}, \quad
v_3 = \begin{pmatrix}
k & l & m \\
0 & 0 & p \\
0 & 1 & n \\
\end{pmatrix}, \quad
v_4 = \begin{pmatrix}
x & y & z \\
0 & 1 & u \\
0 & 0 & t \\
\end{pmatrix}.
$$

We have $v_4^2 = v_4$, and there are the following variants 
of the Jordan form~$J$ of $v_4$:

1${}^\prime$) $J = e_{22}$,

2${}^\prime$) $J = e_{22}+e_{33}$,

3${}^\prime$) $J = e_{11}+e_{22}$.

The case 3${}^\prime$) is conjugate to the second one under~$\Theta_{13}$.

{\sc Case 1${}^\prime$}: $v_4 = e_{22}$.
From $v_2v_4 = v_4v_2 = 0$,
$v_4v_1 = v_1$ and $v_3v_4 = v_3$ we get
$$
v_1 = \begin{pmatrix}
0 & 0 & 0 \\
1 & 0 & d \\
0 & 0 & 0 \\
\end{pmatrix}, \quad
v_2 = \begin{pmatrix}
f & 0 & h \\
0 & 0 & 0 \\
1 & 0 & j \\
\end{pmatrix}, \quad
v_3 = \begin{pmatrix}
0 & l & 0 \\
0 & 0 & 0 \\
0 & 1 & 0 \\
\end{pmatrix}.
$$
By~$v_3v_1 = v_2$, we express
$v_2 = \begin{pmatrix}
l & 0 & dl \\
0 & 0 & 0 \\
1 & 0 & d \\
\end{pmatrix}$.
If $d+l = 0$, then
$\psi$~\eqref{AutU} with parameters
$\beta = \varepsilon = 0$, $\gamma = d$,
$\alpha = \delta = 1$ acts as follows, 
$\psi(e_{21}) = v_1$,
$\psi(e_{31}) = v_2$,
$\psi(e_{32}) = v_3$,
it is~(B6).

If $d+l\neq0$, then
$\psi$~\eqref{AutU} with parameters
$\beta = \varepsilon = 0$, $\gamma = d$,
$\alpha = d+l$, $\delta = 1$ acts as follows, 
$\psi(e_{21}) = v_1$,
$\psi(e_{12}+e_{32}) = v_2/(d+l)$,
$\psi(e_{11}+e_{31}) = v_3/(d+l)$,
we get~(B7) with the help of~$\psi^{-1}$.

{\sc Case 2${}^\prime$}: $v_4 = e_{22}+e_{33}$.
From the equations 
$v_1v_4 = v_2v_4 = 0$, 
$v_3v_4 = v_3$,
$v_4v_i = v_i$, $i=1,2,3$,
we derive that $v_1 = e_{21}$, $v_2 = e_{31}$, and 
$v_3 = \begin{pmatrix}
0 & 0 & 0 \\
0 & 0 & p \\
0 & 1 & n \\
\end{pmatrix}$.
Analogously to the case~2), we get 
either decomposition~(B8), when $p+n^2/4 = 0$,
or~(B9), otherwise.

The group of decompositions (B1)--(B5) and the group of decompositions (B6)--(B9) have non-isomorphic biggest subalgebra~$M$, thus, it is enough to show that decompositions lying in the same group are from different orbits.
Inside the first group, decompositions (B2) and (B4) but not others have 2-dimensional radicals of corresponding subalgebras $S$. Further, the radical of $S$ from (B4) but not from (B2) lies in one-sided annihilator of the whole~$S$.
The rank of the idempotent lying in~$S$ from (B1) equals~1, and the same parameter for (B3) and (B5) equals~2.
Also, there exists an idempotent in $S$ from (B5), which acts as unit on the square of its radical, and there are no such idempotents in~$S$ from the decomposition~(B3).

For the second group, we have $S\cong M_2(\mathbb{C})$ only for (B7).
The semisimple part of~$S$ is 2-dimensional only for~(B9).
Finally, the ranks of the idempotents lying in~$S$ from (B6) and (B8) are not equal.
\hfill $\square$

\section*{Acknowledgments}

Vsevolod Gubarev is supported by RAS Fundamental Research Program, project FWNF-2022-0002.

\noindent Vsevolod Gubarev \\
Sobolev Institute of Mathematics \\
Acad. Koptyug ave. 4, 630090 Novosibirsk, Russia \\
Novosibirsk State University \\
Pirogova str. 2, 630090 Novosibirsk, Russia \\
e-mail: wsewolod89@gmail.com
\end{document}